\documentclass[11pt,reqno]{amsart}
\usepackage{indentfirst, amssymb, amsmath, amsthm, mathrsfs, setspace, indentfirst, enumerate,  mathrsfs, amsmath, amsthm}
\usepackage[bookmarksnumbered, colorlinks, plainpages]{hyperref}
\usepackage{mathrsfs}
\hypersetup{
	colorlinks=true,
	linkcolor=blue,
	citecolor=green,
}
\newtheorem*{theo2A}{Theorem 2.A}
\newtheorem*{theo2B}{Theorem 2.B}
\newtheorem*{theo2C}{Theorem 2.C}
\newtheorem*{theo2D}{Theorem 2.D}
\newtheorem*{theo2E}{Theorem 2.E}
\newtheorem*{theo2F}{Theorem 2.F}
\newtheorem*{theo2G}{Theorem 2.G}
\newtheorem*{theo2H}{Theorem 2.H}
\newtheorem*{theo2I}{Theorem 2.I}

\newtheorem*{ques2A}{Question 2.A}
\newtheorem*{theo3A}{Theorem 3.A}

\newtheorem*{ques3A}{Question 3.A}
\newtheorem*{theo4A}{Theorem 4.A}
\newtheorem*{theo4B}{Theorem 4.B}
\newtheorem*{theo4C}{Theorem 4.C}
\newtheorem*{theo4D}{Theorem 4.D}
\newtheorem*{conj4A}{Conjecture 4.A}

\newtheorem*{exm2A}{Example 2.A}
\newtheorem*{conj2A}{Conjecture 2.A}
\newtheorem{theo}{Theorem}[section]
\newtheorem{lem}{Lemma}[section]
\newtheorem{cor}{Corollary}[section]

\newtheorem{exm}{Example}[section]

\newtheorem{rem}{Remark}[section]
\newtheorem{ques}{Question}[section]

\newcommand{\ol}{\overline}
\newcommand{\be}{\begin{equation}}
\newcommand{\ee}{\end{equation}}
\newcommand{\beas}{\begin{eqnarray*}}
\newcommand{\eeas}{\end{eqnarray*}}
\newcommand{\bea}{\begin{eqnarray}}
\newcommand{\eea}{\end{eqnarray}}

\numberwithin{equation}{section}
\begin{document}
\title[E\MakeLowercase{ntire function related to shift}......]{\Large E\MakeLowercase{ntire function related to shift and difference operator}}
\date{}
\author[S. M\MakeLowercase{ajumder}, L. M\MakeLowercase{ahato and} A. D\MakeLowercase{as}]{S\MakeLowercase{ujoy} M\MakeLowercase{ajumder$^1$$^{*}$},  L\MakeLowercase{ata} M\MakeLowercase{ahato$^2$ and} A\MakeLowercase{shoke} D\MakeLowercase{as}$^3$}

\address{$^1$Department of Mathematics, Raiganj University, Raiganj, West Bengal-733134, India.}
\email{sm05math@gmail.com, sjm@raiganjuniversity.ac.in}

\address{$^2$Department of Mathematics, Gangarampur College, Gangarampur, Dakshin Dinajpur, West Bengal-733124, India.}
\email{lata27math@gmail.com}

\address{$^3$Department of Mathematics, Raiganj University, Raiganj, West Bengal-733134, India.}
\email{ashoke.avik@gmail.com}

\renewcommand{\thefootnote}{}
\footnote{*\emph{Corresponding Author}: Sujoy Majumder.}

\renewcommand{\thefootnote}{\arabic{footnote}}
\setcounter{footnote}{0}

\begin{abstract}In this paper, we investigate the uniqueness problems of entire functions involving difference operators. First, we examine the case where $f-a$, $\Delta_c^n f-b$, and $\Delta_c^{n+2}f-b$ share $0$ CM (counting multiplicities), where $a$ and $b$ are two entire small functions of $f$ such that $\Delta_c^m a \equiv 0$ ($1 \le m \le n$), and $b$ is $c$-periodic. As a consequence, we establish two main results, one of which resolves a conjecture posed by Farissi et al. \cite{FLBA1}. Second, we obtain a uniqueness theorem for an entire function $f(z)$ when $f(z)$, $f(z+c)$, and $\Delta_c f(z)$ share a single value, thereby improving and generalizing the corresponding result of Farissi et al. \cite{FLA1} by removing the growth restriction on the order of $f$. Third, we establish uniqueness results for $f(z)$ when $f(z)$ and $\Delta_c^k f(z)$ share a value, which improve and generalize the findings of Chen et al. \cite{CCL2} and Huang and Zhang \cite{HZ1}. Finally, we provide several examples to demonstrate that our results are best possible in certain senses.
\end{abstract}
\thanks{Typeset by \AmS -\LaTeX}
\maketitle

\section{{\bf Introduction}}

\noindent
The study of entire functions in relation to their shift and difference operators has emerged as an important area of modern complex analysis, particularly within the framework of Nevanlinna value distribution theory. This subject explores how the growth, value distribution, and uniqueness properties of an entire function are connected to those of its shifted and difference expressions. In recent decades, the development of difference analogues of classical differential equations has attracted considerable attention, leading to significant advances in the theory of complex difference equations, functional equations, and discrete dynamical systems.
\medskip

\noindent
Let $c$ be a non-zero finite complex number. Then the forward difference $\Delta^n_c f$ for each positive integer $n$ is defined by 
\begin{align*}
\Delta^1_c f(z)=\Delta_c f(z)=f(z+c)-f(z)
\end{align*}
\begin{align*}
\Delta_c^n f(z)=\Delta_c\left(\Delta_c^{n-1}f(z)\right)=\Delta_c^{n-1} f(z+c)-\Delta_c^{n-1}f(z),\;\;n\geq 2.
\end{align*}
Moreover
\begin{align*}
 \Delta_c^n f(z)=\sum\limits_{j=0}^n(-1)^{n-j}C^{j}_{n}f(z+jc),
\end{align*}
where $C^{j}_{n}$ is a combinatorial number.
\smallskip

\noindent
These operators provide discrete counterparts to differentiation and play a fundamental role in investigating the analytic and growth properties of entire functions. While differentiation describes local variation, the difference operator measures the change in a function over a fixed complex displacement. The relationship between these two operations becomes especially significant when studying entire functions of finite order, exponential type, or more general growth classes.
\medskip

\noindent
A central tool in this investigation is Nevanlinna value distribution theory, which provides a systematic framework for understanding the growth and distribution of the values assumed by entire and meromorphic functions.
One of the principal directions of research concerns the uniqueness of entire functions sharing values with their shifts or differences. In particular, researchers investigate conditions under which the sharing of one or more complex values, either counting multiplicities (CM) or ignoring multiplicities (IM), between \(f\) and \(f(z+c)\), or between \(f\) and \(\Delta_c f\), forces a specific functional relationship. Such problems are closely related to classical uniqueness theorems involving derivatives and have inspired discrete analogues of several fundamental results in Nevanlinna theory. The central question is whether prescribed value-sharing conditions impose sufficient rigidity to determine the form of the entire function. Another important research direction concerns complex difference equations and difference-differential equations.
\smallskip

\noindent
The interaction between shift and difference operators also raises questions concerning periodicity, finite-order growth, and the characterization of exponential functions. For example, if \(f(z+c)=f(z)\), then \(f\) is periodic with period \(c\), whereas the equation \(\Delta_c f(z)=0\) expresses the same periodicity condition.
In addition to these problems, the study of entire functions and their shift and difference operators has important connections with discrete dynamical systems, functional equations, and the theory of meromorphic solutions of nonlinear difference equations. The extension of Nevanlinna theory to difference operators has established a framework for studying these connections through growth estimates, value distribution, and uniqueness principles. For this background, we see \cite{CCL2}-\cite{HK1}, \cite{HZ1}, \cite{HF1}, \cite{LFB1}-\cite{LD1}, \cite{M1}-\cite{QY1}, \cite{WY1}, \cite{HYX}.
\medskip

\noindent
We define an $\varepsilon$-set to be a countable union of discs (see \cite{BL1})
\begin{align*}\displaystyle
 E=\bigcup\limits_{j=1}^{\infty} B(b_j,r_j)
\end{align*}
such that
\begin{align*}\displaystyle
	 \lim\limits_{j\to\infty}|b_j|=\infty\;\;\text{and}\;\;\sum\limits_{j=1}^{\infty}\frac{r_j}{|b_j|}<\infty.
\end{align*}

\noindent
Here $B(a, r)$ denotes the open disc of center $a$ and radius $r$, and $S(a, r)$ will denote the corresponding boundary circle. Note that if $E$ is an $\varepsilon$-set, then the set of $r\geq 1$ for which the circle $S(0, r)$ meets $E$ has a finite logarithmic measure.
\medskip

\noindent
We assume that the reader is familiar with standard notation and main
results of Nevanlinna Theory (see \cite{WKH,YY1}). By $S(r, f)$ we denote any quantity that satisfies the condition $S(r, f) = o(T(r, f))$ as $r\to \infty$ possibly outside of an exceptional set of finite logarithmic measure. A meromorphic function $a$ is said to be a small function of $f$ if $T(r,a)=S(r,f)$ and we denote by $S(f)$ the set of functions which are small compared to $f$.
Moreover, we use notations $\rho(f)$, $\mu(f)$ and $\rho_1(f)$ for the order, lower order and the hyper-order of a meromorphic function $f$ respectively.
Let $f$ be a non-constant meromorphic function in $\mathbb{C}$ and $a\in\mathbb{C}\cup\{\infty\}$. We define
\begin{align*} \delta (a;f)=\liminf _{r \rightarrow \infty}\frac{m(r,a;f)}{T(r,f)} = 1- \limsup _{r \rightarrow \infty }\frac{N(r,a;f)}{T(r,f)}
\end{align*}
The quantities $\delta(a;f)$ is known as the deficiency of $a$. Clearly $\delta(a;f)\in [0,1]$.
\medskip

\noindent
The paper is organized as follows. In Section \ref{sec2}, we will discuss the uniqueness problem of entire function $f$ when $f-a$, $\Delta_c^n f-b$ and $\Delta_c^{n+2}f-b$ share $0$ CM, where $a$ and $b$ are two small entire functions of $f$ such that $\Delta_c^ma\equiv 0\;(1\leq m\leq n)$ and $b$ is $c$-periodic. Also in this section we will try to settle the conjecture posed by Farissi et al. \cite{FLBA1}. In section \ref{sec3}, we will prove a result to show that the condition ``order of $f$ is finite'' is not necessary in a result of Farissi et al. \cite{FLA1}. 
In Section \ref{sec4}, we will improve as well as generalise the results of Chen et al. \cite{CCL2} and  Huang and Zhang \cite{HZ1}.

\section{\bf{Shared value problem of $f$, $\Delta_c^nf$ and $\Delta_c^{n+2}f$}}\label{sec2}
\noindent
In 1977, Rubel and Yang \cite{RY1} established that if an entire function $f$ shares two distinct complex numbers CM (counting multiplicities) with its first derivative $f^{(1)}$, then $f \equiv f^{(1)}$. Subsequently, Jank et al. \cite{JMV1} investigated the case of a single shared value involving a meromorphic function and its derivatives, proving that if a non-constant meromorphic function $f$, along with $f^{(1)}$ and $f^{(2)}$, shares a finite non-zero value CM, then $f \equiv f^{(1)}$.
\medskip

\noindent
In 2001, Li and Yang \cite{LY1} proved the following result.
\begin{theo2A}\cite[Theorem 1]{LY1} Let $f$ be a non-constant entire function, $a$ be a finite non-zero constant and let $n$ be a positive integer. If $f$, $f^{(n)}$ and $f^{(n+1)}$ share the value $a$ CM, then $f\equiv f^{(1)}$.
\end{theo2A}

\noindent
In 2007, Bergweiler and Langley \cite{BL1} proved that: For a  transcendental meromorphic function $f$ of order $\rho(f)<1$, we have
\begin{align*}
	\Delta_c^k f(z)\sim c^kf^{(k)}(z),\quad z\rightarrow \infty\quad \text{in}\;\mathbb{C}\backslash  E,
\end{align*}
where $E$ is an $\varepsilon$-set. 

\noindent
In view of the asymptotic relation between $f^{(k)}$ and $\Delta_c^k f$, it is natural to ask whether a corresponding analogue holds for difference operators, leading to the following question related to Theorem 2.A:

\begin{ques2A} What is the uniqueness result in the case when a finite-order entire function $f$ shares value with $\Delta_c^n f$ and $\Delta_c^{n+1}f$, where $n$ is a positive integer ? 
\end{ques2A}

\noindent
In connection with Question 2.A, Chen et al. \cite{CCL1, CL1} established the following result:
\begin{theo2B}\cite[Theorem 1.1]{CCL1} Let $f$ be a finite order non-constant entire function and let $a\in S(f)$ be a periodic entire function with period $c$ such that $a\not\equiv 0$. If $f$, $\Delta_cf$ and $\Delta_c^2f$ share $a$ CM, then $\Delta_cf\equiv \Delta_c^2f$.
\end{theo2B}

\begin{theo2C}\cite[Theorem 1.1]{CL1} Let $f$ be a finite order non-constant entire function and let $a\in S(f)$ be a periodic entire function with period $c$ such that $a\not\equiv 0$. If $f$, $\Delta_cf$ and $\Delta_c^nf\;(n\geq 2)$ share $a$ CM, then $\Delta_cf\equiv \Delta_c^nf$.
\end{theo2C}

\noindent
It is natural to investigate what happens when $f$, $\Delta_c^n f$, and $\Delta_c^{n+1} f$ ($n \ge 1$) share a value $a$ CM. In this direction, Latreuch et al. \cite{LFB1} established the following result:

\begin{theo2D}\cite[Theorem 1.7]{LFB1} Let $f$ be a finite order non-constant entire function and let $a\in S(f)$ be a periodic entire function with period $c$ such that $a\not\equiv 0$. If $f$, $\Delta_c^nf$ and $\Delta_c^{n+1}f\;(n\geq 1)$ share $a$ CM, then $\Delta_c^{n+1}f\equiv \Delta_c^nf$.
\end{theo2D}

\noindent
In pursuit of a difference analogue of the result by Li and Yang \cite{LY1}, Farissi et al. \cite{FLA1} established the following results:

\begin{theo2E}\cite[Theorem 1.1]{FLA1} Let $f$ be a non-periodic entire function of finite order and let $a\in S(f)$ be a periodic entire function with period $c$ such that $a\not\equiv 0$.  If $f$, $\Delta_cf$ and $\Delta_c^2f$ share $a$ CM, then $\Delta_cf\equiv f$.
\end{theo2E}

\begin{theo2F}\cite[Corollary 1.1]{FLA1} Let $f$ be a non-periodic entire function of finite order and let $a(\not\equiv 0)\in S(f)$ be a periodic entire function with period $c$. If $f$, $\Delta_cf$ and $\Delta_c^3f$ share $a$ CM, then $\Delta_cf\equiv f$.
\end{theo2F}

\noindent
In 2016, Farissi et al. \cite{FLBA1} further improved Theorem 2.E in several directions. We now recall their results:

\begin{theo2G}\cite[Theorem 1.8]{FLBA1} Let $f$ be a finite order non-constant entire function such that $\Delta_c^nf\not\equiv 0 $ and let $a\in S(f)$  be a periodic entire function with period $c$ such that $a\not\equiv 0$. If $f$, $\Delta_c^nf$ and $\Delta_c^{n+1}f\;(n\geq 1)$ share $a$ CM, then $\Delta_cf\equiv f$.
\end{theo2G}

\noindent
The following example shows that the condition $\Delta_c^n f \not\equiv 0$ is necessary in Theorem 2.G:
\begin{exm2A}\cite{FLBA1} Let $f(z)=e^{2\pi i z}+1$ and $a=c=1$. Note that $\Delta_c^n f-a=\Delta_c^{n+1}f-a=-1$. Clearly $f$, $\Delta_c^nf$ and $\Delta_c^{n+1}f$ share $a$ CM, but $\Delta_cf\not\equiv f$.
\end{exm2A}

\begin{theo2H}\cite[Theorem 1.10]{FLBA1} Let $f$ be a finite order non-constant entire function such that $\Delta_c^nf\not\equiv 0 $ and let $a(\not\equiv 0),\;b(\not\equiv 0)\in S(f)$ such that $b$ be a periodic entire function with period $c$ and $\Delta_c^ma\equiv 0(1\leq m\leq n)$. If $f-a$, $\Delta_c^nf-b$ and $\Delta_c^{n+1}f-b$ share $0$ CM, then $\Delta_cf\equiv f+b+\Delta_ca-a.$ 
\end{theo2H}

\noindent
In the case where $b \equiv 0$, Farissi et al. \cite{FLBA1} obtained the following result:
\begin{theo2I}\cite[Theorem 1.17]{FLBA1} Let $f$ be a finite order non-constant entire function such that $\Delta_c^nf\not\equiv 0 $ and let $a\in S(f)$ such that $\Delta_c^ma\equiv 0\;(1\leq m\leq n)$. If $f-a$, $\Delta_c^nf$ and $\Delta_c^{n+1}f$ share $0$ CM, then $\Delta_cf\equiv Cf+\Delta_ca-a,$ where $C$ is a non-zero constant. 
\end{theo2I}

\noindent
At the end of their paper, Farissi et al. \cite{FLBA1} noted that it remains an open question to determine the conditions under which Theorem 2.F holds for entire functions $f(z)$ that share a small function with $\Delta_c^n f(z)$ and $\Delta_c^{n+2} f(z)$ ($n \ge 1$). In this direction, Farissi et al. \cite{FLBA1} posed the following conjecture:

\begin{conj2A}\cite{FLBA1} Let $f$ be a finite order non-constant entire function such that $\Delta_c^nf\not\equiv 0$ and let $a(\not\equiv 0)\in S(f)$ be a periodic entire function with period $c$. If $f$, $\Delta_c^nf$ and $\Delta_c^{n+2}f\;(n\geq 1)$ share $a$ CM, then $\Delta_c f\equiv f$.
\end{conj2A}
\medskip

\noindent
In this paper, we resolve Conjecture 2.A under the additional assumption that $\delta(0;f)>0$. We now state our main results:

\begin{theo}\label{t1} Let $f$ be a finite order non-constant entire function such that $\Delta_c^nf\not\equiv 0$ and let $a, b(\not\equiv 0)\in S(f)$ be entire functions such that $\Delta_c^ma\equiv 0\; (1\leq m\leq n)$ and $b$ is $c$-periodic. If $f-a$, $\Delta_c^nf-b$ and $\Delta_c^{n+2}f-b$ share $0$ CM and $\delta(a-b;f)>0$, then one of the following conclusion holds:
\begin{enumerate}
\item[(i)] $\Delta_cf\equiv f+b+\Delta_ca-a$, if $n$ is odd;
\item[(ii)] $\Delta_c^2f\equiv f+b+\Delta_c^2a-a$, if $n$ is even.
\end{enumerate} 
\end{theo}

\noindent
When $b\equiv 0$, we get the following result.
\begin{theo}\label{t2} Let $f$ be a non-constant entire function such that $\rho(f)<2$ and $\Delta_c^nf\not\equiv 0$ and let $a\in S(f)$ such that $\Delta_c^ma\equiv 0\; (1\leq m\leq n)$. If $f-a$, $\Delta_c^nf$ and $\Delta_c^{n+2}f$ share $0$ CM, then one of the following conclusion holds:
	\begin{enumerate}
		\item[(i)] $\Delta_cf\equiv C(f-a)+\Delta_ca$, if $n$ is odd;
		\item[(ii)] $\Delta_c^2f\equiv C(f-a)+\Delta_c^2a$, if $n$ is even, 
	\end{enumerate} 
where $C$ is a non-zero constant.	
\end{theo}
\noindent
As an immediate consequence of Theorem \ref{t1}, we obtain the following corollary:

\begin{cor}\label{c1} Let $f$ be a finite order non-constant entire function such that $\Delta_c^nf\not\equiv 0$ and $\delta(0;f)>0$ and let $a(\not\equiv 0)\in S(f)$ be a periodic entire function with period $c$. If $f-a$, $\Delta_c^nf-a$ and $\Delta_c^{n+2}f-a$ share $0$ CM, then $\Delta_cf\equiv f$.
\end{cor}

\noindent
The following example demonstrates that condition $(ii)$ in Theorem \ref{t1} cannot be omitted.
\begin{exm} Let $f(z)=e^{\alpha z}-1+z$, $c=1$, $e^{\alpha}=2$, $a(z)=z$ and $b(z)=1$. Note that $\Delta_c^2a\equiv 0$, $\Delta_cf(z)=e^{\alpha z}+1$, $\Delta_c^2f(z)=e^{\alpha z}$ and $\Delta_c^4f(z)=e^{\alpha z}$. Observe that $f-a+b$ has no zeros and so $\delta(a-b;f)=1$. Here
\begin{align*}
\frac{\Delta_c^2f-b}{f-a}=1\;\;\;\text{and}\;\;\;\frac{\Delta_c^4f-b}{f-a}=1.
\end{align*}
Clearly $f-a$, $\Delta_c^2f-b$ and $\Delta_c^4f-b$ share $0$ CM and $f$ satisfies the difference equation $\Delta_c^2f\equiv f+b+\Delta_c^2a-a$. 
\end{exm}

\noindent
As illustrated by the following example, conclusion $(ii)$ of Theorem \ref{t2} is essential and cannot be dropped.

\begin{exm}Let $f(z)=e^{\frac{1}{4}\left(\log 2+\frac{\pi}{2}\iota \right)z}+z$, $c=2$, $n=2$ and $a(z)=z$. Note that 
\begin{align*}
\Delta_cf(z)=&e^{\frac{1}{4}\left(\log 2+\frac{\pi}{2}\iota\right)(z+2)}+(z+2)-e^{\frac{1}{4}\left(\log 2+\frac{\pi}{2}\iota\right)z}-z\\
=&e^{\frac{1}{4}\left(\log 2+\frac{\pi}{2}\iota\right)z}\left(e^{\frac{1}{2}\left(\log 2+\frac{\pi}{2}\iota\right)}-1\right)+2\\
=&e^{\frac{1}{4}\left(\log 2+\frac{\pi}{2}\iota\right)z}\left(\sqrt{2}\left(\frac{1}{\sqrt{2}}+\iota\frac{1}{\sqrt{2}}\right)-1\right)+2\\
=&\iota e^{\frac{1}{4}\left(\log 2+\frac{\pi}{2}\iota\right)z}+2.
\end{align*}
Therefore $\Delta_c^2f(z)=-e^{\frac{1}{4}\left(\log 2+\frac{\pi}{2}\iota\right)z}$ and $\Delta_c^4f(z)=e^{\frac{1}{4}\left(\log 2+\frac{\pi}{2}\iota\right)z}$. It follows that $\Delta_c^2f=-(f-a)$, $\Delta_c^4f=f-a$ and so $f-a$, $\Delta_c^2f$ and $\Delta_c^4f$ share $0$ CM. We see that $f$ satisfies the difference equation
$\Delta_c^2f\equiv C(f-a)+\Delta_c^2a$, where $C=-1$.      
\end{exm}

\noindent
The following example illustrates the existence of entire functions that satisfy conclusion $(i)$ of Theorem \ref{t2}.

\begin{exm} Let $f(z)=e^{\frac{1}{4}\left(\log 2+\frac{\pi}{2}\iota\right)z}+1$, $c=2$, $n=1$ and $a(z)=1$. Note that 
\begin{align*}
\Delta_cf(z)=e^{\frac{1}{4}\left(\log 2+\frac{\pi}{2}\iota\right)(z+2)}-e^{\frac{1}{4}\left(\log 2+\frac{\pi}{2}\iota\right)z}
=&e^{\frac{1}{4}\left(\log 2+\frac{\pi}{2}\iota\right)z}\left(e^{\frac{1}{2}\left(\log 2+\frac{\pi}{2}\iota\right)}-1\right)\\
=&\iota e^{\frac{1}{4}\left(\log 2+\frac{\pi}{2}\iota\right)z}.
\end{align*}
It is clear that $\Delta_cf=\iota (f-a)$, $\Delta_c^3f=-\iota (f-a)$ and so $f-a$, $\Delta_cf$ and $\Delta_c^3f$ share $0$ CM. We see that $f$ satisfies the difference equation $\Delta_cf\equiv C(f-a)+\Delta_ca$, where $C=\iota$.             
\end{exm}

\noindent
\subsection*{\bf{ Auxiliary lemmas}}

\begin{lem}\label{l5}\label{3l2}\cite[Theorem 2.1]{CF1} Let $f$ be a meromorphic function of finite order $\rho$ and let $c$ be a fixed non-zero complex number. Then for any $\varepsilon >0$, we have
	\begin{align*}
		T(r,f(z+c))=T(r,f(z))+O(r^{\rho(f)-1+\varepsilon })+O(\log r).
	\end{align*}
\end{lem}

\noindent
Furthermore, if $f$ is a transcendental meromorphic function of finite order, then
we know
\begin{align*}
	T(r,f(z+c))=T(r,f(z))+S(r,f).
\end{align*}

\begin{lem}\label{3l1} \cite[Theorem 2.1]{HK1} Let $f$ be a finite order transcendentalt meromorphic function, $c$ be a finite complex number. Then
\begin{align*}
m\left(r,\frac{f(z+c)}{f(z)}\right)=S(r,f).
\end{align*}
\end{lem}

\begin{lem}\label{3l5} \cite[Lemma 3]{LW1} Suppose that $h$ is a non-constant meromorphic function satisfying
\begin{align*}
\ol N(r,h)+\ol N(r,1/h)=S(r,h).
\end{align*}
Let $f_1=a_0h^p+a_1h^{p-1}+\cdots+a_p,$ and $f_2=b_0h^q+a_1h^{q-1}+\cdots+a_q,$ be polynomials in h with coefficients $a_0, a_1,\cdots,a_p$, $b_0,b_1,\cdots,b_q$ being small functions of $h$ and $a_0b_0a_p\neq 0.$ If $q\leq p,$ then $m(r,f_2/f_1)=S(r, h).$
\end{lem}

\begin{lem}\label{314}\cite[Lemma 3.3]{BL1} Let $f$ be a transcendental meromorphic function in the plane such that $\rho(f)<1$. Let $h>0$. Then there exists an $\varepsilon $-set $E$ such that
	\begin{align*}
		\frac{f^{(1)}(z+c)}{f(z+c)}\rightarrow 0\;\;\text{and}\;\;\frac{f(z+c)}{f(z)}\rightarrow 1\;\text{as}\; z\rightarrow \infty\;in\;\mathbb{C}\backslash  E
	\end{align*}
	uniformly in c for $|c|\leq h$. Further $E$ may be chosen so that for large $z$ not in $E$ the function $f$ has no zeros or poles in $|\zeta -z|\leq h$.
\end{lem}




\begin{lem}\label{l4} Let $f$ be a finite order transcendental entire function such that $\Delta_c^nf\not\equiv 0$, let $a(\not\equiv 0)$ be an entire small function of $f$ and let $P$ be a polynomial. If $f$ is a solution of the following equation 
\begin{align}\label{l2.1}
\Delta_c^nf-a=e^{P}(f-a),
\end{align}
such that $\delta(0;f)>0$, then $P$ reduces to a constant and $\rho(f)\geq 1$.
\end{lem}

\begin{proof} Re-writing \eqref{l2.1}, we get 
\begin{align}\label{l2.1a} 
\frac{1}{f}=\frac{1}{a(1-e^P)}\frac{\Delta_{c}^n f}{f}-\frac{e^P}{a(1-e^P)}.
\end{align}
We know that
\begin{align*}
	\frac{\Delta_{c}^n f(z)}{f(z)}=\sum\limits_{j=0}^n(-1)^{n-j}C^{j}_{n}\frac{f(z+jc)}{f(z)}
\end{align*}
and so by Lemma \ref{3l1}, we obtain
\begin{align}\label{l2.2}
	m\left(r,\frac{\Delta_{c}^n f}{f}\right)=S(r,f).
\end{align}
We now consider the following two cases:
\medskip

\noindent
{\bf Case 1.} Let $\rho(f)<1$. We know that $\rho(\Delta_c^nf)\leq \rho(f)$. Then $\rho(\Delta_c^nf-a)<1$ and $\rho(f-a)<1$ and so from \eqref{l2.1}, we have $\rho(e^P)<1$, which shows that $P$ is a constant. Let $e^P=d$. 
\smallskip

\noindent
Suppose $d=1$. Then from (\ref{l2.1}), we get $\Delta_c^nf=f$ and  so
\begin{align}\label{l2.2a}
	\sum\limits_{j=0}^n(-1)^{n-j}C^{j}_{n}\frac{f(z+jc)}{f(z)}=1.
\end{align}
Using Lemma \ref{314} to \eqref{l2.2a}, we get $(1-1)^n=1$, which is impossible.
\smallskip

\noindent 
Next suppose $d\neq 1$. Consequently \eqref{l2.1a}, gives
\begin{align}\label{l2.3} 
\frac{1}{f}=\frac{1}{a(1-d)}\frac{\Delta_{c}^n f}{f}-\frac{d}{a(1-a)}.
\end{align}
 Therefore using \eqref{l2.2} to \eqref{l2.3}, we find that
$m(r,1/f)=S(r,f)$. Hence $T(r,f)=N(r,0;f)+S(r,f)$ and so $\delta(0;f)=0$, which is a contradiction. 
\medskip

\noindent
{\bf Case 2.} Let $\rho(f)\geq 1$. If possible, suppose $P$ is non-constant. 
With the help of Lemma \ref{3l5},  we deduce that 
\begin{align}\label{l2.4}
m\left(r,\frac{1}{1-e^P}\right)=S(r,e^P)\;\;\text{and}\;\;m\left(r,\frac{e^P}{1-e^P}\right)=S(r,e^P).
\end{align}
Applying Lemma \ref{l5} to (\ref{l2.1}), we get $T(r,e^P)=O(T(r,f))$. This shows that $S(r,e^P)$ can be replaced by $S(r,f)$. Now using \eqref{l2.2} and \eqref{l2.4} to \eqref{l2.3}, we find that
$m(r,1/f)=S(r,f)$ and so $\delta(0;f)=0$, which is a contradiction. Hence $P$ is a constant. 
\end{proof}


\noindent
An example is given to show the necessity of the condition $\delta(0;f)>0$ in Lemma \ref{l4}.
\begin{exm} 
The function $f(z)=1+2e^{2\pi \iota z}-e^{\pi \iota z}$ satisfies the difference equation
$\Delta_1 f(z)-1=e^{-\pi \iota z}(f(z)-1)$. Note that $T(r,f)=2T(r,e^{\pi \iota z})+O(1)$. Suppose $f(z)=2(e^{\pi \iota z}-\alpha_1)(e^{\pi \iota z}-\alpha_2)$.
Clearly $N(r,\alpha_1;e^{\pi \iota z})\sim T(r,e^{\pi \iota z})$ and $N(r,\alpha_2;e^{\pi \iota z})\sim T(r,e^{\pi \iota z})$. But, we see that
\begin{align*} \delta(0;f)=&1-\lim\limits_{r\rightarrow \infty}\sup\frac{N(r,0;f)}{T(r,f)}\\
=&1-\lim\limits_{r\rightarrow \infty}\sup\frac{N(r,\alpha_1;e^{\pi \iota z})+N(r,\alpha_2;e^{\pi \iota z})}{2T(r,e^{\pi \iota z})+O(1)}\\
=&1-\lim\limits_{r\rightarrow \infty}\sup\frac{2T(r,e^{\pi \iota z})}{2T(r,e^{\pi \iota z})+O(1)}=0.
\end{align*}
\end{exm}

\begin{lem}\label{l1} \cite[Corollary 2.6]{CF1} Let $f$ be a meromorphic function of finite order $\rho$ and let $\eta_1,\;\eta_2\in\mathbb{C}$ such that $\eta_1\neq \eta_2$. Then for any $\varepsilon>0$, we have
	\begin{align*} 
		m\left(r,\frac{f(z+\eta_1)}{f(z+\eta_2)}\right)=O(r^{\rho-1+\varepsilon}).
	\end{align*}
\end{lem}

\medskip
\begin{proof}[{\bf Proof of Theorem \ref{t1}}]
By the given conditions, we know that $f-a$, $\Delta_c^nf-b$ and $\Delta_c^{n+2}f-b$ share $0$ CM, where $\Delta_c^ma=0$ and $\delta(a-b;f)>0$.
If possible, suppose that $f$ is a non-constant polynomial of degree $n_1$. We know that a small function of a non-constant polynomial is a constant. Therefore $a$ and $b$ are constants. Clearly $\Delta_c^n f$ is a polynomial of degree $n_1-n$.
Therefore $f-a$ and $\Delta_c^nf-b$ cannot have the same zeros with the same multiplicities. So, we get a contradiction.
Hence we assume that $f$ is a transcendental entire function.

\noindent
Let $g=f+b-a$. Then $g-b$, $\Delta_c^ng-b$ and $\Delta_c^{n+2}g-b$ share $0$ CM and $\delta(0;g)>0$. Since $\rho(f)=\rho(g)<+\infty$, by Hadamard's factorization theorem there exist two polynomials $P$ and $Q$ such that
\begin{align}\label{t3.1} 
\frac{\Delta_c^n g-b}{g-b}=e^{P}\;\;\text{and}\;\;\frac{\Delta_c^{n+2}g-b}{g-b}=e^{Q}.
\end{align}
Applying Lemma \ref{l4} to (\ref{t3.1}), we find that $P$ and $Q$ are both constants.
Let $e^P=d_1$ and $e^Q=d_2$. Then from (\ref{t3.1}), we have   
\begin{align}\label{t3.4}
\Delta_c^n g=d_1g+(1-d_1)b
\end{align}
and
\begin{align}\label{t3.5}
\Delta_c^{n+2} g=d_2g+(1-d_2)b.
\end{align}
Since $b$ is a periodic function with period $c$, from (\ref{t3.4}), we get $\Delta_c^{n+2}g=d_1\Delta_c^2 g$ and so from (\ref{t3.5}), we have
\begin{align}\label{t3.8} 
\Delta_c^2g=\frac{d_2}{d_1}g+\frac{1-d_2}{d_1}b.
\end{align}

\noindent
The analysis now splits into two distinct cases:
\medskip

\noindent
{\bf Case 1.} Let $n$ be odd. 

\noindent
Suppose $d=\frac{d_2}{d_1}$. Then by induction and using (\ref{t3.8}) repeatedly, we obtain the following 
\begin{align}\label{t3.9} 
\Delta_c^ng=d^{\frac{n-1}{2}}\Delta_cg,\;\text{if}\;\;n\geq 3
\end{align}
and
\begin{align}\label{t3.10}
\Delta_c^{n+2}g=d^{\frac{n+1}{2}}\Delta_cg,\;\text{if}\;\;n\geq 1.
\end{align}
We note that (\ref{t3.9}) is also true for $n=1$. Using (\ref{t3.4})-(\ref{t3.5}) and (\ref{t3.9})-(\ref{t3.10}), we respectively have
\begin{align}\label{t3.12} 
d^{\frac{n-1}{2}}\Delta_cg-d_1g=(1-d_1)b
\end{align}
and
\begin{align}\label{t3.13} 
d^{\frac{n+1}{2}}\Delta_cg-d_2g=(1-d_2)b.
\end{align}
Now from (\ref{t3.12}) and (\ref{t3.13}), we get $\big(d^{\frac{n+1}{2}}-d^{\frac{n-1}{2}}\big)= 0$ and so $d=1$, i.e., $d_1=d_2$. Therefore from (\ref{t3.4}) and (\ref{t3.9}), we get
\begin{align}\label{t3.15} 
\Delta_c g=d_1g+(1-d_1)b.
\end{align}
Again from (\ref{t3.8}) and (\ref{t3.15}), we get $\Delta_c^3g=\Delta_cg$ and $\Delta_c^3g=d_1^2\Delta_cg$ respectively and so $d_1^2=1$, i.e., $d=\pm 1$. If $d=-1$, then from (\ref{t3.15}), we get
\begin{align}\label{t3.16} 
\Delta_c g=-g-2b,
\end{align}
which implies $g(z+c)=-2b(z)$. Thus using Lemma \ref{l5}, we get 
\begin{align*}
T(r,f)=T(r,f(z+c))+S(r,f)=T(r,g(z+c))+S(r,f)=S(r,f),
\end{align*}
which is a contradiction. If $d=1$, then from (\ref{t3.15}), we get $\Delta_c g=g$, i.e., $\Delta_cf=f+b+\Delta_ca-a$.
\smallskip

\noindent
{\bf Case 2.} Let $n$ be even. 

\noindent
First suppose $n=2$. Since $b$ is a periodic function with period $c$, from (\ref{t3.4}), we get $\Delta_c^{4}g=d_1\Delta_c^2 g$ and so from (\ref{t3.5}), we obtain
\begin{align}\label{t3.17} 
d_1\Delta_c^2g=d_2g+(1-d_2)a.
\end{align}
Now (\ref{t3.4}) and (\ref{t3.17}) yield
\begin{align}\label{t3.18} 
(1-d_1)\Delta_c^2g=(d_1-d_2)g.
\end{align}
We claim that $d_1=1$. If not, suppose $d_1\neq 1$. Then from (\ref{t3.4}) and (\ref{t3.18}), we get
\begin{align}\label{t3.19}
\left(d_1-\frac{d_1-d_2}{1-d_1}\right)g=(d_1-1)b,
\end{align}
which shows that $d_1\neq \frac{d_1-d_2}{1-d_1}$. Then from (\ref{t3.19}), we get $T(r,g)=T(r,f)+S(r,f)=S(r,f)$, which is impossible. Hence $d_1=1$ and so from (\ref{t3.18}), we have $d_1=d_2=1$. Consequently from (\ref{t3.17}), we get $\Delta_c^2g=g$, i.e., $\Delta_c^2f=f+b+\Delta_c^2a-a$.

\noindent
Next suppose $n\geq 4$. Then by induction and using (\ref{t3.8}) repeatedly, we obtain the following 
\begin{align}\label{t3.20} 
\Delta_c^ng=d^{\frac{n}{2}}g
\end{align}
and
\begin{align}\label{t3.21}
\Delta_c^{n+2}g=d^{\frac{n+2}{2}}g.
\end{align}
Using (\ref{t3.4})-(\ref{t3.5}) and (\ref{t3.20})-(\ref{t3.21}), we get respectively
\begin{align}\label{t3.22} 
\left(d^{\frac{n}{2}}-d_1\right)g=(1-d_1)b\;\;
\text{and}\;\; 
\left(d^{\frac{n+2}{2}}-d_2\right)g=(1-d_2)b.
\end{align}
Clearly from (\ref{t3.22}), we get $(1-d_1)b\equiv 0$ and $(1-d_2)b\equiv 0$. Since $b\not\equiv 0$, we have $d_1=d_2=1$. Therefore from (\ref{t3.8}), we get $\Delta_c^2g=g$, i.e., $\Delta_c^2f=f+b+\Delta_c^2a-a$. 
\end{proof}

\begin{proof}[{\bf Proof of Theorem \ref{t2}}]
 By the given conditions, we know that $f-a$, $\Delta_c^nf$ and $\Delta_c^{n+2}f$ share $0$ CM and $\Delta_c^ma=0$. Let $g=f-a$. Then $g$, $\Delta_c^ng$ and $\Delta_c^{n+2}g$ share $0$ CM. Since $\rho(f)=\rho(g)<2$, by Hadamard's factorization theorem there exist two polynomials $P$ and $Q$ such that
\begin{align}\label{s3.1} 
\frac{\Delta_c^n g}{g}=e^{P}\;\;\text{and}\;\;\frac{\Delta_c^{n+2}g}{g}=e^{Q}.
\end{align}
Clearly \eqref{s3.1} gives
\begin{align}\label{s3.2}
e^{P(z)}=\sum\limits_{j=0}^n(-1)^{n-j}C^{j}_{n}\frac{g(z+ic)}{g(z)}.
\end{align}
Let us choose $\varepsilon>0$ such that $\rho(g)+\varepsilon<2$. Now using Lemma \ref{l1} to \eqref{s3.2}, we find that
\begin{align*}
T\left(r,e^P\right)=m\left(r,e^P\right)=O\left(r^{\rho(g)-1+\varepsilon}\right)+O(1),
\end{align*}
which demands that $P$ reduces to a constant. Similarly, we can prove that $Q$ is also a constant. Let $e^P=d_1$ and $e^Q=d_2$. Then from \eqref{s3.1}, we get
\begin{align}\label{s3.4}
\Delta_c^n g=d_1g
\end{align}
and
\begin{align}\label{s3.5}
\Delta_c^{n+2} g=d_2g.
\end{align}

\noindent
Let $d=\frac{d_2}{d_1}$. Then from (\ref{s3.4}), we get $\Delta_c^{n+2}g=d_1\Delta_c^2 g$ and so from (\ref{s3.5}), we arrive at (\ref{t3.8}), where $b\equiv 0$.
\medskip

\noindent
First, assume that $n$ is odd. Now using (\ref{t3.10}) and (\ref{s3.5}), we get $\Delta_cg=Cg$, i.e., $\Delta_cf=C(f-a)+\Delta_ca$, where $C=d_2/d^{\frac{n+1}{2}}$.
\medskip

\noindent
Next, assume that $n$ is even. Since $\Delta_c^{n+2}g=d_1\Delta_c^2 g$, from (\ref{s3.5}), we get $\Delta_c^2g=dg$, i.e., $\Delta_c^2f=d(f-a)+\Delta_c^2a$.
\end{proof}

\section{\bf{Shared value problem of $f(z)$, $f(z+c)$ and $\Delta_cf(z)$}}\label{sec3}
\noindent
In 2016, Farissi et al. \cite{FLA1} raised the following question:
\begin{ques3A} What will happen if $f(z)$, $f(z + c)$ and $\Delta_cf(z)$ share a finite value $a$ CM ?
\end{ques3A}

\noindent
In \cite{FLA1}, Farissi et al. completely answered Question 3.A by establishing the following result:
 
\begin{theo3A}\cite[Theorem 1.4]{FLA1} Let $f$ be a finite order transcendental entire function such that $\Delta_cf\not\equiv 0$. Then $f(z)$, $f(z+c)$ and $\Delta_cf$ can not share any finite value $a\neq 0$ CM. Furthermore, if $a=0$, then $f$ must be of the form $f(z)=h(z)e^{\frac{d}{c}z}$, where $d\neq 0$ and $h(z)$ is periodic entire function of period $c$.
\end{theo3A}

\noindent
In this paper, we extend and refine Theorem 3.A along two main directions:
\begin{enumerate}
\item[(1)] We remove the condition that $f$ is of finite order.
\item[(2)] We consider shared small functions instead of shared values.
\end{enumerate}
Below, we present our result:

\begin{theo} \label{3t1} Let $f$ be a transcendental entire function such that $\Delta_cf\not\equiv 0$ and let $a$ be a $c$-periodic small entire function of $f$. Then $f(z)-a(z)$, $f(z+c)-a(z)$ and $\Delta_cf(z)-a(z)$ can not share $0$ CM if $a\not\equiv 0$. Furthermore, if $a\equiv 0$, then $f(z)$ must be of the form $f(z)=h(z)e^{\frac{d}{c}z}$ where $d\neq 0$ is a finite complex number and $h(z)$ is a $c$-periodic entire function.
\end{theo}

\noindent
\subsection*{\bf {Auxiliary lemma}}
\begin{lem}\label{l2} \cite[Theorem 1.64]{YY1} Let $f_1,\ldots, f_n$ be non-constant meromorphic functions and $f_{n+1},\ldots, f_{n+m}$ be meromorphic functions such that $f_k\not\equiv 0\;\;(k=n+1,\ldots, n+m)$ and 
\begin{align*}
\sum_{i=1}^{n+m} f_i\equiv A,	
\end{align*}
a non-zero constant. If there exists a subset $I\subseteq \mathbb{R}^+$ satisfying $\text{mes}\;I=+\infty$ such that
\begin{align*}
\sum\limits_{i=1}^{n+m} N(r,0;f_i)+(n+m-1)\sum\limits_{i=1, i\neq j}^{n+m} \ol N(r,f_i)<(\lambda+o(1))T(r,f_j)\;\;(r\to \infty, r\in I, j=1,2,\ldots,n),
\end{align*}
where $\lambda<1$, then there exist $t_i\in\{0,1\}\;(i=1,2,\ldots,m)$ such that $\sum_{i=1}^mt_if_{n+i}\equiv A$.
\end{lem} 

\begin{proof}[{\bf Proof of Theorem \ref{3t1}}]
By the given conditions, we know that $f(z)-a(z)$, $f(z+c)-a(z)$ and $\Delta_cf(z)-a(z)$ share $0$ CM. Then there exist entire functions $P$ and $Q$ such that
\begin{align}\label{t2.1} 
\frac{f(z+c)-a(z)}{f(z)-a(z)}=e^{P(z)}\;\;\text{and}\;\;\frac{\Delta_cf(z)-a(z)}{f(z)-a(z)}=e^{Q(z)}.
\end{align}

\noindent
We consider the following two sub-cases:
\medskip

\noindent
{\bf Case 1.} Let $a\equiv 0$. Then from (\ref{t2.1}), we get
\begin{align}\label{t2.3} 
e^{P}-e^{Q}=1,
\end{align}
which shows that either both $P$ and $Q$ are constant or both are non-constant. If $P$ and $Q$ are both non-constant, then (\ref{t2.3}) implies that $1$ is a Picard exceptional value of $e^P$. Consequently, $e^P$ would possess two Picard exceptional values ($0$ and $1$), which is impossible for an entire function. Hence, both $P$ and $Q$ must be constant. Setting $e^P = e^d$, we obtain from (\ref{t2.1}) that
\begin{align}\label{t2.4} 
f(z+c)=e^df(z).
\end{align}
Let $f$ and $g$ be any two solutions of (\ref{t2.4}). If we take $h=\frac{f}{g}$, then (\ref{t2.4}) gives $h(z+c)=h(z)$ and so
$h$ is a $c$-periodic function. On the other hand, it is clear that $g(z)=e^{\frac{d}{c}z}$ is a solution of (\ref{t2.4}). Hence the solution of (\ref{t2.4}) is of the form $f(z)=h(z)e^{\frac{d}{c}z}$, where $h$ is $c$-periodic entire function.
\medskip

\noindent
{\bf Case 2.} Let $a\not\equiv 0$. Set $g=f-a$. Then from (\ref{t2.1}), we get 
\begin{align}\label{t2.5} 
g(z+c)=g(z)e^{P(z)}\;\;\text{and}\;\;g(z+c)-g(z)-a(z)=g(z)e^{Q(z)},
\end{align}
which imply that
\begin{align}\label{t2.7}
g(z)=\frac{a(z)}{e^{P(z)}-e^{Q(z)}-1}
\end{align}
and so
\begin{align}\label{t2.8} 
g(z+c)=\frac{a(z)}{e^{P(z+c)}-e^{Q(z+c)}-1}.
\end{align}
Obviously (\ref{t2.5}) and (\ref{t2.7}) yield
\begin{align}\label{t2.9} 
g(z+c)=\frac{a(z)e^{P(z)}}{e^{P(z)}-e^{Q(z)}-1}.
\end{align}
Also from (\ref{t2.8}) and (\ref{t2.9}), we get the following 
\begin{align}\label{t2.10} 
2e^{P(z)}-e^{Q(z)}-e^{P(z)+P(z+c)}+e^{P(z)+Q(z+c)}=1.
\end{align}
Set 
\begin{align*}
h_1(z)=2e^{P(z)},\;\; h_2(z)=-e^{Q(z)},\;\; h_3(z)=-e^{P(z)+P(z+c)}\;\;\text{and}\;\;h_4(z)=e^{P(z)+Q(z+c)}.
\end{align*}
Note that 
\begin{align*}
\sum\limits_{i=1}^{4} N(r,0;h_i)+2\sum\limits_{i=1,i\neq j}^{4} \ol N(r,h_i)<(\lambda+o(1))T(r,h_j)\;\;(r\to \infty, r\in I, j=1,2),
\end{align*}
$\lambda<1$. Then by Lemma \ref{l2}, there exist $t_i\in\{0,1\}, i=1,2$ such that $t_1h_3+t_2h_4=1$, i.e.,
\begin{align}\label{t2.11} 
-t_1e^{P(z)+P(z+c)}+t_2e^{P(z)+Q(z+c)}=1.
\end{align} 
We can see that $(t_1,t_2)\neq (0,0)$. 
\medskip

\noindent
First suppose $t_1=0$. Then from (\ref{t2.11}), we get $e^{P(z)+Q(z+c)}=1$ and so from (\ref{t2.10}), we have $e^{P(z)}(2-e^{P(z+c)})=e^{Q(z)}$, which shows that $N(r,2;e^{P(z+c)})=S(r,e^{P(z+c)})$. Now using the second fundamental theorem, we get 
$T(r,e^{P(z+c)})=S(r,e^{P(z+c)}),$ which is a contradiction.
\medskip

\noindent
Next suppose $t_2=0$. Then from (\ref{t2.11}), we get $e^{P(z)+P(z+c)}=-1$ and so from (\ref{t2.10}), we have $e^{P(z)}(2+e^{Q(z+c)})=e^{Q(z)}$, which shows that $N(r,-2;e^{Q(z+c)})=S(r,e^{Q(z+c)})$. Again using the second fundamental theorem, we get $T(r,e^{Q(z+c)})=S(r,e^{Q(z+c)})$,
which is impossible. 
\medskip

\noindent
Hence $t_1=1$ and $t_2=1$. Now from (\ref{t2.10}) and (\ref{t2.11}), we get $2e^P=e^Q$ and so from (\ref{t2.7}), we have
\begin{align}\label{t2.12} 
f=\frac{ae^P}{e^P-1}.
\end{align}
Clearly $T(r,f)+S(r,f)=T(r, e^P)+S(r,e^P)$, which gives $S(r, f)=S(r,e^P)$. Since $f$ is an entire function, from (\ref{t2.12}), we deduce that $N(r,1;e^P)=N(r,0;a)$ and so by the second fundamental theorem, we get $T(r,e^{P}) \leq N(r,0,a)+S(r,e^P)=S(r,e^P)$, which is a contradiction. 
\end{proof}

\section{\bf{Shared value problem of $f$ and $\Delta_c^k f$}}\label{sec4}

\noindent
What is the possible relation between an entire function $f$ and its first derivative $f^{(1)}$ if they share a finite value CM? Addressing this question, Br\"{u}ck \cite{RB1} formulated the following conjecture:

\begin{conj4A} Let $f$ be a non-constant entire function such that $\rho_1(f)\not\in\mathbb{N}\cup\{\infty\}$.
If $f$ and $f^{(1)}$ share one finite value $a$ CM, then $f^{(1)}-a=c(f-a)$, where $c$ is a non-zero constant. 
\end{conj4A}

\noindent
In 1996, Br\"{u}ck \cite{RB1} established a partial solution to Conjecture 3.A. Although Conjecture 4.A remains unsettled in its full generality, it has sparked a rich line of research concerning the uniqueness problems of entire functions and their derivatives. Notably, Yang and Zhang \cite{YZ1} demonstrated that a variant of Conjecture 4.A holds under appropriate power conditions, proving the following theorem:

\begin{theo4A}\cite[Theorem 4.4]{YZ1} Let $f$ be a non-constant entire function and $n \geq 7$ be an integer. If $f^{n}$ and  $(f^{n})^{(1)}$ share $1$ CM, then $f^{n}\equiv (f^{n})^{(1)}$ and $f(z)=ce^{z/n}$, where $c$ is a non-zero constant. 
\end{theo4A}

\noindent
In 2010, Zhang and Yang \cite{ZY1} proved the following results, significantly extending Theorem 4.A:

\begin{theo4B}\cite[Theorem 2.1]{ZY1} Let $f$ be a non-constant entire function and let $k$ and $n$ be positive integers such that $n \geq k + 1$. If $f^{n}$ and $(f^{n})^{(k)}$ share $1$ CM, then $f^{n}\equiv (f^{n})^{(k)}$ and $f(z)=ce^{\lambda z/n}$, where $c$ and $\lambda$ are non-zero constant such that $\lambda^{k} = 1$. 
\end{theo4B} 
\smallskip

\noindent
We know that $\Delta_c^k f\sim c^kf^{(k)}$ as $z\rightarrow \infty$ in $\mathbb{C}\backslash E$ holds for a transcendental meromorphic function $f$ of order less that $1$, where $E$ is an $\varepsilon$-set. 
\smallskip

\noindent
A natural question therefore arises: what can be said about the uniqueness of a finite-order entire function $f$ when $f^n$ and $\Delta_c^k f^n$ share a value, where $n$ is a positive integer? Addressing this problem, Chen et al. \cite{CCL2} established the following result:

\begin{theo4C}\cite[Theorem 1.5]{CCL2} Let $f$ be a non-constant entire function of finite order, $n\geq 6$ be an integer and let $F=f^n$. If $F$ and $\Delta_cF$ share $1$ CM, then $F\equiv \Delta_cF$.
\end{theo4C}

\noindent
Next, we recall the following result due to Huang and Zhang \cite{HZ1}:

\begin{theo4D}\cite[Theorem 3.4]{HZ1} Let $f$ be a transcendental entire function of finite order, let $a(\not\equiv 0)$ be a small function of $f$ and $c$ is a non-zero finite complex number be such that $\Delta_c f\not\equiv 0$. If $f$ and $\Delta_cf$ share $a$ CM and $\delta(0;f)>0$, then $f\equiv \Delta_cf$.
\end{theo4D}

\noindent
We observe from Theorem 4.D that if the condition $\delta(0;f)>0$ is assumed in Theorem 4.C, then the assertion holds irrespective of the value of $n$. Motivated by Theorem 4.C, it is natural to pose the following question:

\noindent
\begin{ques} Can the lower bound of $n$ in Theorem 4.C be reduced? \end{ques}

\noindent
In this paper, we answer Question 4.1 under the hypothesis that $\delta(0;F)>0$. Furthermore, we extend and improve Theorems 4.C and 4.D along the following directions:
\begin{enumerate}
\item[(1)] We replace value sharing with the sharing of a small function in Theorem 4.C.
\item[(2)] We replace $\Delta_c F$ ($\text{resp.}\ \Delta_c f$) in Theorem 4.C ($\text{resp.}\ \text{Theorem 4.D}$) by $\Delta_c^k F$ ($\text{resp.}\ \Delta_c^k f$).
\end{enumerate}
In this section, we present our main results:

\begin{theo}\label{4t1} Let $f$ be a finite order non-constant entire function, $k$, $n$ be positive integers such that $n\geq k(k+1)+1$ and let $a(\not\equiv 0)$ be an entire small function of $f$ and $c$ be a non-zero finite complex number such that $\Delta_c^k F\not\equiv 0$, where $F=f^n$. Suppose $F-a$ and $\Delta_c^kF-a$ share $0$ CM. If $\delta(0;F)>0$,
then $\rho(f)\geq 1$ and $F\equiv \Delta_c^kF$. In particular when $k=1$, then $f(z)=h(z)e^{\frac{\log 2}{n c}z}$, where $h$ is $c$-periodic entire function.
\end{theo}

\begin{rem} Setting $k=1$ in Theorem \ref{4t1} yields $n \geq 3$; hence, Theorem \ref{4t1} improves Theorem 4.C under the assumption that $\delta(0;F)>0$.
\end{rem}

\noindent
The following example shows that the conditions $n \geq k(k+1)+1$ and $\delta(0;F)>0$ are sharp in Theorem \ref{4t1}.

\begin{exm} Let $f(z)=e^{\frac{z}{2}}+e^{-\frac{z}{2}}$, $c=\pi \iota$, $k=1$, $n=2$ and $a(z)=4/3$. Note that $F(z)=e^z+e^{-z}+2$, where $F(z)=f^2(z)$  and so $\Delta_cF(z)=-2(e^z+e^{-z})$. Obviously $\delta(0;F)=0$. Here we also observe that 
$\Delta_cF-4/3=-2\left(F-4/3\right).$
Clearly $F$ and $\Delta_cF$ share $4/3$ CM, but $F\not\equiv \Delta_cF$.
\end{exm}

\noindent
In the next result, we replace $\Delta_c f$ with the higher-order difference operator $\Delta_c^k f$ in Theorem 4.D:

\begin{theo}\label{4t2} Let $f$ be a finite order non-constant entire function, $a(\not\equiv 0)$ be an entire small function of $f$ and let $c$ be a non-zero constant such that $\Delta_c^k f\not\equiv 0$, where $k$ is a positive integer. If $f-a$ and $\Delta_c^kf-a$ share $0$ CM and $\delta(0;f)>0$, then $\rho(f)\geq 1$ and $f\equiv \Delta_c^kf$. In particular, if $k=1$, then $f(z)=h(z)e^{\frac{\log 2}{c}z}$, where $h$ is $c$-periodic entire function.
\end{theo} 

\noindent
Clearly, Theorem \ref{4t2} is an improvement of Theorem 4.D. The following example demonstrates the necessity of the condition $\delta(0;f) > 0$ in Theorem \ref{4t2}:

\begin{exm} Let $f(z)=2e^{z\log 3}+\frac{1}{2}$. Note that $\rho(f)=1$, $\Delta_1f-1=2(f-1)$ and $\delta(0;f)=0$. Clearly $f$ and $\Delta_1f$ share $1$ CM, but $f\not\equiv \Delta_1f$.
\end{exm}

\noindent
\subsection*{\bf{Auxiliary lemma}}

\begin{lem}\label{3l3} \cite[Theorem 1]{NT1} Let $f_0, f_1,\ldots, f_p\;(p\geq 1)$ be $p+1$ non-constant entire functions and let $a_0, a_1\ldots,a_p$ be $p+1$ non-zero meromorphic functions such that $T(r,a_j)=o(T(r,f_j))$ for $j=0,1,\ldots,p$.
Then if the following functional equation 
\begin{align*}
 a_0f_0^{k_0}+a_{1}f_{1}^{k_1}+\cdots+a_pf_p^{k_p}=1,
 \end{align*}
holds for some integers $k_0,k_1,\cdots, k_p(\geq 1)$, it must be 
\begin{align*}
\sum_{j=0}^p\frac{1}{k_j}\geq \frac{1}{p}.
\end{align*}
\end{lem}

\medskip
\begin{proof}[{\bf Proof of Theorem \ref{4t1}}]
Note that $F-a$ and $\Delta_c^kF-a$ share $0$ CM and $\delta(0;F)>0$, where $F=f^n$. It is easy to deduce that $f$ is a transcendental entire function.
By Hadamard's factorization theorem there exists a polynomial $P$ such that
\begin{align}\label{e3t.1} 
\frac{\Delta_c^k F-a}{F-a}=e^{P}.
\end{align}
Now applying Lemma \ref{l4} to (\ref{e3t.1}), we find that $\rho(f)\geq 1$ and $P$ is a constant.
 Let $e^P=d$. Then from (\ref{e3t.1}), we get
\begin{align}\label{e3t.2} 
\Delta_c^k F-dF=(1-d)a.
\end{align}

\noindent
First we suppose $d=1$. Clearly from (\ref{e3t.2}), we get $F\equiv \Delta_c^kF$. 
In particular, if $k=1$, then we have  $F(z+c)=2F(z)$, i.e., $f(z+c)=e^{\frac{\log 2}{n}}f(z)$. Now proceeding in the same way as done in the proof of Theorem \ref{3t1}, we get $f(z)=h(z)e^{\frac{\log 2}{n c}z}$, where $h$ is $c$-periodic entire function.
\smallskip

\noindent
Next we suppose $d\neq 1$. Again from (\ref{e3t.2}), we get
\begin{align}\label{e3t.3} 
&F(z+kc)+(-1)^{1}C^{1}_{k}F(z+(k-1)c)+\cdots+(-1)^{k-1}C^{k-1}_{k}F(z+c)\nonumber\\
&+((-1)^{k}-d)F(z)=(1-d)a(z).
\end{align}

\noindent
Let $d\neq (-1)^k$. Then from (\ref{e3t.3}), we get
\begin{align}\label{e3t.4} 
a_kf_k^n+a_{k-1}f_{k-1}^n+\cdots+a_1f_1^n+a_0f_0^n=1,
\end{align}
where $a_0=\frac{(-1)^{k}-d}{(1-d)a}$, $a_i=\frac{(-1)^{k-i}C^{i}_{k}}{(1-d)a}$ for $i=1,2,\ldots, k$ and $f_i(z)=f(z+ic)$ for $i=0,1,\ldots, k$. 
\smallskip

\noindent
Since $T(r,a)=o(T(r,f))$, we have $T(r,a_i)=o(T(r,f))$ for $i=0,1,\ldots, k$. Also by Lemma \ref{3l2}, we have $T(r,f_i(z))=T(r,f(z+ic))=T(r,f(z))+S(r,f(z))$ for $i=0,1,\ldots, k$ and so $T(r,a_i)=o(T(r,f_i))$, $i=0,1,\ldots,k$.
Using Lemma \ref{3l3} to (\ref{e3t.4}), we get 
\begin{align*}
\sum\limits_{i=0}^k\frac{1}{n}\geq \frac{1}{k},\;\; \text{i.e.,}\;\; n\leq (k+1)k,
\end{align*}
which contradicts the fact that $n\geq k(k+1)+1$.

\noindent
Let $d=(-1)^k$. Then from (\ref{e3t.4}), we get 
\begin{align}\label{e3t.5} a_kf_k^n+a_{k-1}f_{k-1}^n+\cdots+a_1f_1^n=1,\end{align}
where $T(r,a_i)=o(T(r,f_i))$, $i=1,\ldots,k$.
Applying Lemma \ref{3l3} to (\ref{e3t.5}), we have
\begin{align*}
\sum\limits_{i=1}^k\frac{1}{n}\geq \frac{1}{k-1},\;\; \text{i.e.,}\;\; n\leq (k-1)k,
\end{align*}
which is impossible.
\end{proof}

\begin{proof}[{\bf Proof of Theorem \ref{4t2}}] By the given condition, we know that $f-a$ and $\Delta_c^kf-a$ share $0$ CM and $\delta(0;f)>0$. Obviously $f$ is a transcendental entire function. Now by Hadamard Factorization theorem there exists a polynomial $P$ such that 
\begin{align}\label{s1} 
\frac{\Delta_c^kf-a}{f-a}=e^P.
\end{align}
Applying Lemma \ref{l4} to (\ref{s1}), we find that $\rho(f)\geq 1$ and $P$ reduces to a constant.
 If $e^P=d\neq 1$, then using  \eqref{l2.2} to (\ref{l2.3}), we obtain $m(r,1/f)=S(r,f)$ and so $\delta(0;f)=0$, which is impossible. Hence $e^P=1$ and so from (\ref{s1}), we get $\Delta_c^kf=f$. 
 In particular, if $k=1$, then we have $f(z)=h(z)e^{\frac{\log 2}{c}z}$, where $h$ is $c$-periodic entire function.
\end{proof}

\vspace{1.3mm}

\noindent \textbf {Conflict of interest:} The authors declare that there are no conflicts of interest regarding the publication of this paper.\vspace{1.5mm}

\noindent{\bf Funding:} There is no funding received from any organizations for this research work.
\vspace{1.5mm}

\noindent \textbf {Data availability statement:}  Data sharing is not applicable to this article as no database were generated or analyzed during the current study.
\vspace{1.5mm}


\end{document}